\documentclass{article}
\usepackage{graphicx} 
\usepackage{amsfonts}
\usepackage{amsmath}
\newtheorem{theorem}{Theorem}
\newtheorem{definition}{Definition}
\newtheorem{lemma}
{Lemma}
\usepackage{cite}
\title{Hitting Time Distributions of Random Walks on Finite Graphs}
\author{Anuraag Kumar}
\date{April 2025}

\begin{document}

\maketitle
\begin{abstract}
We investigate the hitting times of random walks on graphs, where a hitting time is defined as the number of steps required for a random walker to move from one node to another. While much of the existing literature focuses on calculating or bounding expected hitting times, this approach is insufficient, as hitting time distributions often exhibit high variance. To address this gap, we analyze both the full distributions and variances of hitting times. Using general Markov chain techniques, as well as Fourier and spectral methods, we derive formulas and recurrence relations for computing these distributions.
\end{abstract}
\section*{Introduction}
Random walks are a fundamental topic in probability theory with wide-ranging applications across physics, engineering, mathematics, and computer science. There are tens of thousands of papers that explore them in depth. Random walks can be used to model a variety of real-world phenomena, including fluid dynamics, stock price movements, genetic drift, animal foraging behavior, and search algorithms.

The central focus of this paper is the \textit{hitting time}, which refers to the number of steps required for a random walk to move from one vertex of a graph to another. In formal terms, for a graph \( G = (V,E) \), we define the hitting time as
\[
\tau_{i,j} = \inf\{t \mid X_0 = i, X_t = j\}
\]
where \( X_t \) denotes the position of the random walker at time \( t \), and \( i \) and \( j \) are two distinct vertices in the graph. 

One effective way to explore hitting times is through simulation. For example, consider the cycle graph with 10 nodes. The following diagram shows the graph:

\begin{center}
    \includegraphics[scale = 0.3]{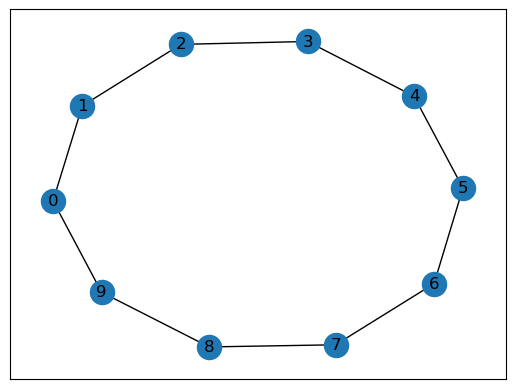}
\end{center}

Suppose we are interested in the random variable \( \tau_{0,5} \), representing the number of steps needed to go from node 0 to node 5. After running 10,000 trials, we obtain the following frequency graph:

\begin{center}
    \includegraphics[scale = 0.3]{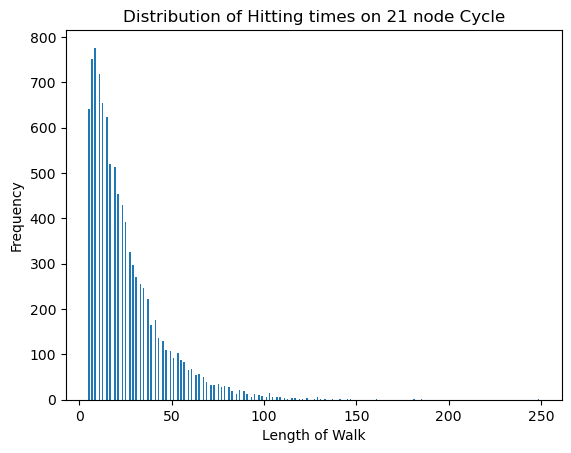}
\end{center}

What stands out here is the \textit{large variance} in the distribution. The sample mean is 25.0306, and the sample variance is 410, indicating that the variance is relatively large compared to the mean.

We can observe a similar behavior in the hypercube graph with 8 nodes, which is shown below:

\begin{center}
    \includegraphics[scale = 0.3]{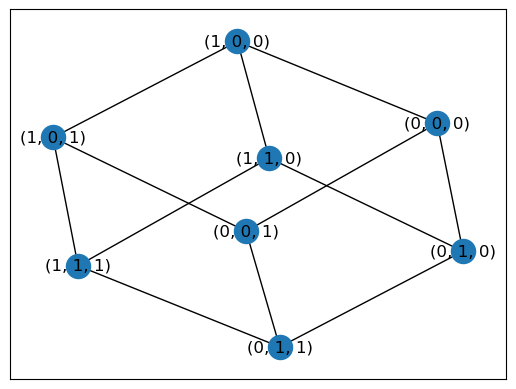}
\end{center}

Now, suppose we are interested in \( \tau_{(0,0,0), (1,1,1)} \), the number of steps needed to move from node \( (0,0,0) \) to node \( (1,1,1) \). After running 10,000 trials, we obtain the following frequency graph:

\begin{center}
    \includegraphics[scale = 0.3]{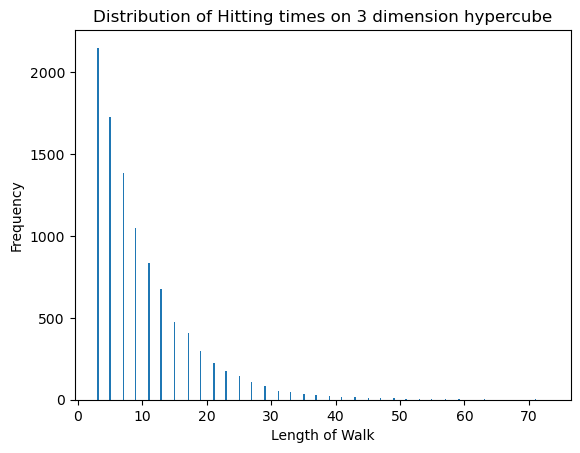}
\end{center}

The sample mean in this case is 10, and the sample variance is 63.

This observation leads to the central motivation for this paper. In the literature, the primary focus in the study of hitting times has been on \textit{expected hitting times}, denoted \( E[\tau_{i,j}] \)\cite{Laneve2023hittingtimesgeneral}\cite{tanaka2023average} \cite{xia2024meanshittingtimesrandom}\cite{zhang2023formulashittingtimescover}. However, as demonstrated above, the expected value alone is often a poor predictor of how the random walk behaves. In this paper, we will examine two related quantities: \( P(\tau_{i,j} = n) \), the probability mass function of hitting times, and \( \text{Var}(\tau_{i,j}) \), the variance of the hitting time. Despite their importance, these quantities receive less attention due to computational challenges. Our approach will involve general computations, with a subsequent focus on vertex-transitive graphs.
 
\section*{Hitting Times on General Graphs}
\subsection*{Distributions}
Let's try to find the distribution of $P(\tau_{i,j} = n)$ on our graph $G$ with Markov matrix $A$. One relationship becomes clear. 
\[P(\tau_{i,j} = n) = \sum_{k \ne j}P(\tau_{i,k} = 1)P(\tau_{k,j} = n-1)\]
The above formula calculates the probability of reaching an adjacent node to our ending node in $n-1$ steps and then making a step from that adjacent node to the end. By setting $k \ne j$, we make sure that we aren't adding the probability that we arrive to ending node $j$ one move early.
The above is true, as each step of a random walk is independent. 
Let us fix an ending node $j$, and then let us define a vector
\[P_n = \begin{bmatrix} 
P(\tau_{1,j} = n) \\
P(\tau_{2,j} = n) \\
\vdots \\
P(\tau_{|V| - 1,j} = n)
\end{bmatrix}\]
the nodes $1,2...|V|-1$ represent some arbitrary numbering of the nodes of the graph once $j$ is removed. Let $Q$ be the matrix such that $Q_{ik} = P(\tau_{ik} = 1)$ such that $i,k \ne j$. We then have
\[P_n = QP_{n-1}\]
as the recursion above is simply matrix multiplication.
So then by induction, we have the following theorem
\begin{theorem}[The Hitting Time Recurrence Formula]
    \[P_n = Q^{n-1}P_1\]
\end{theorem}
One might notice that $Q$ is simply the Markov matrix of our graph but with removed jth row and jth column. Therefore, here we have a solid way of calculating distributions. As we are taking arbitrary powers of a matrix, it often comes down to diagonalizing $Q$. This is often hard to do by hand, but a computer can help. There are a couple things to note about the eigenvalues of such a matrix. \begin{enumerate}
    \item $Q$ is a substochastic matrix as the sum of every row is less than or equal to 1. This implies that $|\lambda| < 1$. Where $\lambda$ is an eigenvalue of $Q$.
    \item $Q$ is the adjacency graph of a subgraph of $G$. Which implies by the Interlacing theorem of Spectral Graph Theory that all the eigenvalues of $Q$ are embedded between the eigenvalues of $A$. 
\end{enumerate}
Consider this graph 
\begin{center}
    \includegraphics[scale = 0.3]{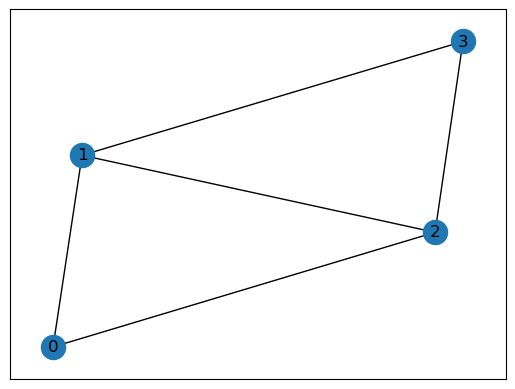}
\end{center}. We will use every step on a random walk is independent and assuming that walking across any edge is equally likely. Setting $j = 0$ and then $Q$ is 
\[
\begin{bmatrix}
0 & \frac{1}{3} & \frac{1}{3} \\
\frac{1}{3}  & 0 & \frac{1}{3} \\
\frac{1}{2} & \frac{1}{2} & 0 \\
\end{bmatrix}
\]
and  $P_1$ is 
\[
\begin{bmatrix}
\frac{1}{3}\\
\frac{1}{3} \\
 0 \\
\end{bmatrix}
\]. 
Then 
\[Q^{n-1}P_1 = \begin{bmatrix}
P(\tau_{1,0} = n) \\
P(\tau_{2,0} = n) \\
P(\tau_{3,0} = n) 
\end{bmatrix} = \begin{bmatrix}
\frac{-\left(-1\right)^n*2^{2n}*3^n*\sqrt{13}+\left(2\sqrt{13}+14\right)^n*\sqrt{13}}{13*\left(\sqrt{13}+1\right)^n2^n3^n} \\
\frac{-\left(-1\right)^n*2^{\left(2*n\right)}*3^n*\sqrt{13}+\left(2\sqrt{13}+14\right)^n\sqrt{13}}{13*\left(\sqrt{13}+1\right)^n2^n3^n} \\
\frac{13*\left(-1\right)^n*2^{2n}3^n+13\left(2\sqrt{13}+14\right)^n+\left(-1\right)^n2^{2n}3^n\sqrt{13}-\left(2\sqrt{13}+14\right)^n\sqrt{13}}{26*\left(\sqrt{13}+1\right)^n2^n3^n}
\end{bmatrix}\]
We can see that even for relatively simple looking graphs, the distributions can be very complicated and often intractable to compute by hand with larger graphs. Later, we will restrict the graphs we will work with to make sure this process is simpler. 
\subsection*{Characteristic Function}
One technique we consider is that of the characteristic function 
\begin{definition}[Characteristic function]
The characteristic function of a hitting time distribution $\phi_{\tau_{i,j}}(t)$ is defined as follows. 
\[\phi_{\tau_{i,j}}(t)=E[e^{it\tau_{i,j}}]=\sum_{n = -\infty}^\infty e^{int}P(\tau_{i,j} = n) =\sum_{n = 0}^\infty e^{int}P(\tau_{i,j} = n) \]
\end{definition}
This can help us compute moments. 
\begin{theorem}[Fast Formulas for the Moments of Hitting Times] The first two moments of the hitting time distribution can be expressed as follows. \cite{li2019hitting}
\[E[\tau_{*,j}] = (I - Q)^{-1}\mathbf{1} = \sum_{n = 0}^{\infty} Q^n \mathbf{1}\]
\[E[\tau_{*,j}^2] = 2\sum_{n=0}^\infty nQ^n\mathbf{1}\]
$\tau_{*,j}$ is a vector of random variable such that the ith coordinate is the random variable is $\tau_{i,j}$. $Q$ is the matrix such that the jth column and row removed corresponding with the fact that j is the fixed ending point.
\end{theorem}
The Proof is in \textbf{Appendix A}  and also in \cite{li2019hitting}

When it comes to actually computing these quantities. The above calculations exist for approximate values found by a computer. In some cases, it is feasible to use the above formulae to find the distributions for general classes of graphs. The cases below are generally easier to compute. Later in this paper, we will introduce machinery to tackle harder cases. We can summarize these easy cases through the following theorem
\subsection*{Hitting Times on Standard Graphs}

We illustrate hitting time distributions for several canonical graphs.

\subsubsection*{Complete Graph ($K_k$)}
For a complete graph with $k$ nodes, a simple random walk moves uniformly to any other node at each step. The hitting time $\tau_{i,j}$ from node $i$ to node $j$ follows a geometric distribution with success probability $1/(k-1)$:
\[
P(\tau_{i,j} = n) = \left(\frac{k-2}{k-1}\right)^{n-1} \frac{1}{k-1}.
\]
Consequently, the expectation and variance are
\[
E[\tau_{i,j}] = k-1, \quad Var[\tau_{i,j}] = (k-1)(k-2).
\]
Intuitively, each step either moves toward the target (success) or not (failure), producing the geometric pattern.

\subsubsection*{Complete Bipartite Graph ($K_{k_1,k_2}$)}
A complete bipartite graph consists of two disjoint sets $A$ and $B$ of sizes $k_1$ and $k_2$, with edges only between sets. Due to this structure, hitting times between sets occur in \textit{odd steps}, while hitting times within a set occur in \textit{even steps}. Specifically:
\begin{itemize}
    \item For $i \in A$, $j \in B$:
    \[
    P(\tau_{i,j} = 2n-1) = \left(1 - \frac{1}{k_2}\right)^{n-1} \frac{1}{k_2}, \quad 
    E[\tau_{i,j}] = 2k_2 - 1, \quad Var[\tau_{i,j}] = 4 k_2(k_2-1).
    \]
    \item For $i, j \in A, i \neq j$:
    \[
    P(\tau_{i,j} = 2n) = \left(1 - \frac{1}{k_2}\right)^{n-1} \frac{1}{k_2}, \quad 
    E[\tau_{i,j}] = 2k_2, \quad Var[\tau_{i,j}] = 4 k_2(k_2-1).
    \]
\end{itemize}
These results highlight how bipartite structure induces alternating-step patterns in hitting times.

\subsubsection*{Cycle Graph ($C_k$)}
For a $k$-cycle, each node is connected to two neighbors. The hitting time distribution can be expressed using powers of the transition matrix. Once the jth row and column are removed. The remaining $Q$ matrix is a tridiagonal Toeplitz matrix which has an easy diagonalization. This lets us prove the following theorem.
\begin{theorem}[Hitting Times on a k-Cycle]
For a cycle with k-nodes, we have that 
\[P(\tau_{i,j} = n) = \frac{1}{k}\sum_{j = 0}^{k-1} \cos(\frac{j\pi}{k})^{n-1}(\sin(\frac{j\pi}{k}) + \sin(\frac{j(k-1)\pi}{k}))\sin(\frac{ij\pi}{k}) \] 
\[E[\tau_{i,j}] = |i - j|(k - |i - j|) \]
where $|i - j|$ is the positive difference between i and j. This result is also proven in \cite{zhang2023formulashittingtimescover}
\end{theorem}
Detailed proofs are provided in \textbf{Appendix A}.

\subsubsection*{Path Graphs and Cycle Symmetry}

An important class of graphs with tractable hitting times is the \textit{path graph} $P_k$, which consists of $k$ nodes connected sequentially. For concreteness, we take one endpoint as the absorbing state (i.e., the target node).  

\begin{center}
    \includegraphics[scale=0.3]{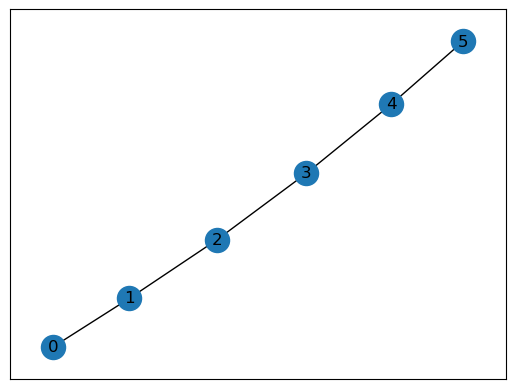}
\end{center}

This choice is without loss of generality: for an interior target, walks that avoid the target can be treated as independent sub-paths where the ending point partitions the graph, so endpoint hitting times capture the general behavior.

We can relate path graphs to cycles using a symmetry argument. Consider a cycle $C_{2k}$ with $2k$ nodes. By reflecting nodes across the line passing through the target, nodes equidistant from the target are paired. Under this mapping, each step toward or away from the target is preserved, and the resulting walk is statistically identical to a random walk on a path of $k$ nodes with an absorbing endpoint.  

\begin{center}
    \includegraphics[scale=0.3]{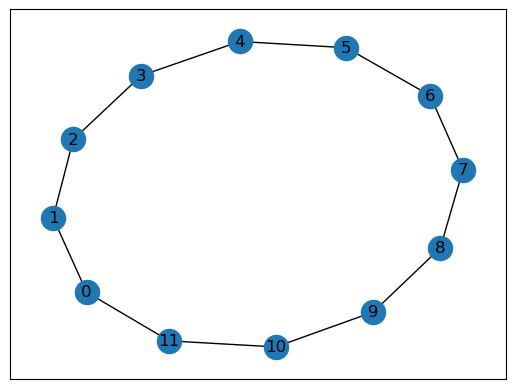}
\end{center}

Thus, the hitting time distribution for a path of $k$ nodes corresponds exactly to that of a cycle of $2k$ nodes under this reflection. This approach allows us to leverage results for cycles to analyze hitting times on paths efficiently. Therefore, Theorem 3 can be used very easily to analyze hitting times on paths.

\section*{Cayley Graphs}
Before we talk about Cayley Graphs. Let us first define a group and related morphisms.
\begin{definition}[Group]
A \textbf{Group} \(G\) is a set equipped with a binary operation \(\cdot\) such that:
\begin{itemize}
    \item (Closure) For all \(g_1, g_2 \in G\), \(g_1 \cdot g_2 \in G\).
    \item (Associativity) For all \(g_1, g_2, g_3 \in G\), \((g_1 \cdot g_2) \cdot g_3 = g_1 \cdot (g_2 \cdot g_3)\).
    \item (Identity) There exists an element \(e \in G\) such that \(e \cdot g = g \cdot e = g\) for all \(g \in G\).
    \item (Inverse) For each \(g \in G\), there exists \(g^{-1} \in G\) such that \(g \cdot g^{-1} = g^{-1} \cdot g = e\).
\end{itemize}
\end{definition}
\begin{definition}[Homomorphism, Isomorphism, and Automorphism]
Let \(G\) and \(H\) be groups.  
\begin{itemize}
    \item A \textbf{homomorphism} is a function \(\phi: G \to H\) that preserves the group operation:
    \[
    \phi(g_1 g_2) = \phi(g_1) \phi(g_2) \quad \forall g_1, g_2 \in G.
    \]
    \item An \textbf{isomorphism} is a bijective homomorphism. If such a map exists, \(G\) and \(H\) are \textbf{isomorphic}, denoted \(G \cong H\). Isomorphic groups are structurally identical, differing only in element labels.
    \item An \textbf{automorphism} is an isomorphism from a group to itself, \(\phi: G \to G\). Automorphisms capture the symmetries of a group and form the \textbf{automorphism group} \(\mathrm{Aut}(G)\) under composition.
\end{itemize}
\end{definition}
Let's suppose we have a random walk on a group $G$. We can imagine each element in our group $g \in G$ being associated with a node in a graph. Let's suppose that two nodes $g,h$ are connected if for some $c \in C \subset G$, $gc = h$. We often call $C$ the connection set. The resultant graph is called a Cayley graph. Many times we want to consider a graph connected only by its generators. This changes for different kinds of graphs. If we want to perform a random walk on such a graph, we assign a probability from moving from $g\to gc$. Let us call this $p(c)$. Generally, $C$ is symmetric. This means that if $x \in C$ then $x^{-1} \in C$. If $p(x) = p(x^{-1})$, then we have a symmetric random walk.

We analyze the Cayley Graphs

In the context of Cayley graphs, the cycle graph on \( k \) nodes can be seen as the Cayley graph of the finite cyclic group \( \mathbb{Z}_k \) with generating set \( \{\pm1\} \). Extending this idea, we now consider the group \( \mathbb{Z}_p^2 \), where \( p \) is an odd prime. A standard set of generators for its Cayley graph is:
\[
\{(\pm1, 0), (0, \pm1)\},
\]
which gives rise to the 2D torus graph — a grid with wrap-around edges.

However, an alternative set of generators is:
\[
\{(\pm1, \pm1)\},
\]
which corresponds to diagonal steps in the lattice. These generators also produce a valid Cayley graph, albeit with a rotated geometry. However, why would we do this? This shift let's us consider random walks that are independent in both coordinates, something that is not true if we are moving along the standard generators of $Z_p^2$. This happens as the probability of moving vertically is inevitably tied to the probability of moving to the side as a walker must pick one or the other. By considering a diagonal random walk, with the above generators, we introduce independence in the walks in either of the coordinates. Importantly, these generator sets are related via a linear automorphism:
\[
\phi(a,b) = \left(\frac{a + b}{2}, \frac{a - b}{2}\right),
\]
which is an automorphism of \( \mathbb{Z}_p^2 \) since \( 2 \) is invertible in \( \mathbb{Z}_p \) (as \( p \) is odd). This transformation maps the standard coordinate basis to the new one spanned by the diagonal generators.

Let \( c_n(i) \) denote the probability that a walk on a \( p \)-cycle starting at \( 0 \) hits \( i \) for the first time at step \( n \). If we define \( \phi^{-1}(a - c, b - d) = (a', b') \), then \( (a', b') \) represents the displacement between \( (a,b) \) and \( (c,d) \) expressed in terms of the diagonal generator basis. Under the assumption of independence in each coordinate (which holds due to the structure of the walk), the hitting time distribution on \( \mathbb{Z}_p^2 \) can be expressed as:
\[
P(\tau_{(a,b),(c,d)} = n) = \sum_{i=0}^n c_i(a') \cdot c_{n-i}(b').
\]
That is, the distribution of the hitting time is a convolution of the 1D hitting time distributions along the transformed coordinates. The reason that $Z_p^2$ was singled out instead of the direct products of other cyclic groups is that the $\phi$ is an automorphism on $Z_p^2$ and not in other such groups. However, if one was to draw those graphs generated by these diagonal generators, the above framework would easily be able to find the distribution for those graphs as well.
\section*{Fourier View}
For the second part of this paper, we will consider random walks on groups. To analyze these walks effectively, we need to introduce some algebraic machinery. We now introduce the Fourier transform on Finite Grounds. 
\subsection*{Fourier Transform on Finite Groups}
We first have a set of definitions to introduce.
\begin{definition}[Group Representation]
A \textbf{representation} of a group \(G\) is a homomorphism \(\rho: G \to \mathrm{GL}(V)\), where \(V\) is a complex vector space and \(\mathrm{GL}(V)\) is the group of invertible linear transformations on \(V\). That is, for all \(g_1, g_2 \in G\):
\[
\rho(g_1 g_2) = \rho(g_1)\rho(g_2).
\]
\end{definition}

\begin{definition}[Irreducible Representation]
A representation \(\rho: G \to \mathrm{GL}(V)\) is \textbf{irreducible} if there is no nontrivial subspace \(W \subset V\) (i.e., \(W \ne \{0\}\) and \(W \ne V\)) such that \(\rho(g)(w) \in W\) for all \(g \in G\) and \(w \in W\). Intuitively, an irreducible representation cannot be decomposed into smaller representations; it is a “building block” for all representations of \(G\).
\end{definition}

\begin{definition}[Trivial Representation]
The \textbf{trivial representation} of a group \(G\) is the one-dimensional representation \(\rho_1\) defined by
\[
\rho_1(g) = 1 \quad \text{for all } g \in G.
\]
\end{definition}

\begin{definition}[Fourier Transform on a Group]
Let \(f: G \to \mathbb{C}\) be a complex-valued function on a finite group \(G\). The \textbf{Fourier transform} of \(f\) at a representation \(\rho\) of \(G\) is defined as
\[
\widehat{f}(\rho) = \sum_{g \in G} f(g) \rho(g).
\]
\end{definition}

\begin{definition}[Abelian Group and Irreducible Representations]
A group \(G\) is \textbf{Abelian} if \(g_1 g_2 = g_2 g_1\) for all \(g_1, g_2 \in G\).  
Important facts about finite Abelian groups:
\begin{itemize}
    \item All irreducible representations are one-dimensional.
    \item The number of irreducible representations equals the order \(|G|\) of the group.
\end{itemize}
\end{definition}
\cite{Fulton:1991}

\begin{definition}[Inverse Fourier Transform]
Let \(G\) be a finite group, and let \(R\) be the set of irreducible representations of \(G\). For a function \(f: G \to \mathbb{C}\), the \textbf{inverse Fourier transform} is given by
\[
f(g) = \frac{1}{|G|} \sum_{\rho \in R} d_\rho \, \mathrm{Tr}\big(\rho(g^{-1}) \widehat{f}(\rho)\big),
\]
where \(d_\rho\) is the dimension of the representation \(\rho\).  

For Abelian groups, all irreducible representations are one-dimensional (\(d_\rho = 1\)), so this simplifies to
\[
f(g) = \frac{1}{|G|} \sum_{\rho \in R} \rho(g^{-1}) \widehat{f}(\rho).
\]
\end{definition}

\begin{lemma}[Plancherel's Theorem]
Let \(f, g: G \to \mathbb{C}\) be functions on a finite group \(G\) with irreducible representations \(R\). Then
\[
\sum_{a \in G} f(a) g(a^{-1}) = \frac{1}{|G|} \sum_{\rho \in R} d_\rho \, \mathrm{Tr} \big( \widehat{f}(\rho) \widehat{g}(\rho) \big).
\]
For Abelian groups, this reduces to
\[
\sum_{a \in G} f(a) g(a^{-1}) = \frac{1}{|G|} \sum_{\rho \in R} \widehat{f}(\rho) \widehat{g}(\rho).
\]
\end{lemma}

\begin{lemma}[Sum over Non-Trivial Representations]
For a non-trivial irreducible representation \(\rho \ne \rho_1\) of a finite group \(G\):
\[
\sum_{g \in G} \rho(g) = 0.
\]
This result is proven in \cite{zhang2023formulashittingtimescover}.
\end{lemma}
The most important reason we introduce this machinery into this paper is the fact that Fourier Transform turns convolutions into products. If for all $a \in G$.
\[h(a) = \sum_{s\in G} f(a)g(as^{-1}) \]
\[\widehat{h}(\rho) = \widehat{f}(\rho)\widehat{g}(\rho) \]

In our case of hitting time distributions, we have that for $x,y \in G$. $p^*(x,y)$ is the probability of moving from x to y. For random walks on groups, we are assuming a time-independent increment distribution. Therefore, we can define a new function $p(g) = p^*(x,xg) = P(\tau_{x,xg} = 1)$. 

\title{Second moment and variance of hitting times on Abelian Cayley graphs}
\author{ }
\date{ }
\maketitle

\section*{Setup}
Let \(G\) be a finite abelian group (written multiplicatively) and let \(p:G\to[0,1]\) be the one-step transition law of a random walk on the Cayley graph of \(G\). Denote by \(\tau_{e,g}\) the first-passage time from the identity \(e\) to \(g\in G\). Define
\[
h(g)=\mathbb{E}[\tau_{e,g}],\qquad q(g)=\mathbb{E}[\tau_{e,g}^2],
\]
and the return-second-moment
\[
q^*=\mathbb{E}[(\tau_e^+)^2].
\]
Write \(R=\{\rho_1,\dots,\rho_{|G|}\}\) for the (one-dimensional) irreducible characters of \(G\), with \(\rho_1\) the trivial character. For each character \(\rho\) set
\[
\widehat p(\rho)=\sum_{s\in G} p(s)\rho(s).
\]

\begin{lemma}[Expected Hitting Time in an Abelian Cayley Graph]
    In \cite{zhang2023formulashittingtimescover}, Zhang proves that if $G$ is the Cayley Graph of an Abelian Group $\mathcal{G}$ then for $i,j \in \mathcal{G}$we have that 
    \[E[\tau_{i,j}] = h(ij^{-1})= \sum_{m = 2}^{|\mathcal{G}|} \frac{1- \rho_m(ij^{-1})}{1 - \widehat{p}(\rho_m)}\]
    where $|\mathcal{G}|$ is the number of elements of the group. $\rho_mj$ is the mjth irreducible representation of $G$. As $G$ is finite and abelian, it has as many irreducible representations as it has elements. The one representation missing from the sum is the trivial one $\rho_1$. 
\end{lemma}

\section*{Recurrence for the second moment}
Conditioning on the first step yields
\[
q(g)=\sum_{s\in G} p(s)\,\mathbb{E}\big[(\tau_{e,gs^{-1}}+1)^2\big]
=\sum_{s\in G} p(s) q(gs^{-1})+2\sum_{s\in G} p(s) h(gs^{-1})+1,
\]
with boundary condition \(q(e)=0\). Introduce the adjustment function
\[
k(g)=\begin{cases}1,&g\neq e,\\[4pt]1-q^*,&g=e,\end{cases}
\]
so the recurrence becomes
\[
q(g)=\sum_{s\in G} p(s) q(gs^{-1}) + 2\sum_{s\in G} p(s) h(gs^{-1}) + k(g).
\]

Taking Fourier transforms (over characters \(\rho_j\)) and using that irreducible characters are one-dimensional for abelian \(G\), for each \(j\ge2\) we obtain
\[
\widehat q(\rho_j)\big(1-\widehat p(\rho_j)\big)
= -2\widehat p(\rho_j)\widehat h(\rho_j)-\widehat k(\rho_j).
\]
For the nontrivial characters \(j\ge2\) one has \(\widehat k(\rho_j)=-q^*\). Substituting the spectral expression for \(\widehat h(\rho_j)\) (consistent with \eqref{eq:h-form}) and solving for \(\widehat q(\rho_j)\) yields
\[
\widehat q(\rho_j)
= \frac{2|G|\,\widehat p(\rho_j)}{\big(1-\widehat p(\rho_j)\big)^2}
-\frac{q^*}{1-\widehat p(\rho_j)},\qquad j\ge2.
\]

Invert the Fourier transform (and use \(\widehat q(\rho_1)=-\sum_{j=2}^{|G|}\widehat q(\rho_j)\) which follows from \(q(e)=0\)) to get a compact expression.

\begin{theorem}[Variance of Hitting Time Distributions]
\[
q(g)=\frac{1}{|G|}\sum_{j=2}^{|G|}\left(
\frac{2|G|\,\widehat p(\rho_j)}{(1-\widehat p(\rho_j))^2}
-\frac{q^*}{1-\widehat p(\rho_j)}
\right)\bigl(1-\rho_j(g^{-1})\bigr).
\]

Since \(h(g)\) is given by Lemma 3, the variance is
\[
\mathrm{Var}[\tau_{e,g}]
=q(g)-h(g)^2
=\frac{1}{|G|}\sum_{j=2}^{|G|}\left(
\frac{2|G|\,\widehat p(\rho_j)}{(1-\widehat p(\rho_j))^2}
-\frac{q^*}{1-\widehat p(\rho_j)}
\right)\bigl(1-\rho_j(g^{-1})\bigr)
-\left(\sum_{j=2}^{|G|}\frac{1-\rho_j(g^{-1})}{|G|\bigl(1-\widehat p(\rho_j)\bigr)}\right)^{\!2}.
\]
It must be noted that $q^* =-\frac{1}{|G|}+\frac{2}{|G|^2} + Z_{jj},$ where $Z = (I-P+11^T\frac{1}{|G|})^{-1}$ by \cite{KemenySnell1960}. 
\end{theorem}

\subsection*{Distributions}
For the ease of writing, we will define that $P(\tau_{g,e} = n) = m_n(g)$. It follows that $p(g) = m_1(g)$
It follows that we have 
\[m_n(g) = \sum_{s \in G} p^*(g,gs^{-1})m_{n-1}(gs^{-1})\]
and 
\[m_n(e) = 0\]
for $n \geq 1$.
So, we can design a new function $c_n$
\[c_n(g) = \begin{cases} 
      0 & g \ne e \\
      \sum_{s \in G} m_{n-1}(s)m_1(s) & g = e \\ 
   \end{cases}
\]
\[m_n(g) = -c_n(g) + \sum_{s \in G} m_1(s^{-1})m_{n-1}(gs^{-1})\]
For the purposes of this write-up, we will assume that our random walk is symmetric. 
\[m_n(g) = -c_n(g) + \sum_{s \in G} m_1(s)m_{n-1}(gs^{-1})\]
\[\widehat{m_n}(\rho_j) = -I\sum_{s \in G} m_{n-1}(s)p(s)+ \widehat{p}(\rho_j)\widehat{m_{n-1}}(\rho_j)\]
\[\widehat{m_n}(\rho_j) = -I\sum_{s \in G} m_{n-1}(s)m_1(s^{-1})+ \widehat{m_1}(\rho_j)\widehat{m_{n-1}}(\rho_j)\]
By Lemma 1 and the fact that Abelain groups have exactly the same number of irreducible representations and group elements, we can say the following.
\[\widehat{m_n}(\rho_j) = -\frac{1}{k}\sum_{a = 0}^{k-1} \widehat{m_{n-1}}(\rho_a)\widehat{m_1}(\rho_a)+ \widehat{m_1}(\rho_j)\widehat{m_{n-1}}(\rho_j)\]

and so then we have that
\[\widehat{m_n}(\rho_j) = -\frac{1}{k}\sum_{a = 0}^{k-1} \widehat{m_{n-1}}(\rho_a)\widehat{m_1}(\rho_a)+ \widehat{m_1}(\rho_j)\widehat{m_{n-1}}(\rho_j)\]
So then we have a recurrence relation such that.
\[A\begin{bmatrix}
    \widehat{m_{n-1}(\rho_0)} \\
    \widehat{m_{n-1}(\rho_1)} \\
    ... \\
    \widehat{m_{n-1}(\rho_{k-1})} 
    \end{bmatrix}= 
    \begin{bmatrix}
    \widehat{m_{n}(\rho_0)} \\
    \widehat{m_{n}(\rho_1)} \\
    ... \\
    \widehat{m_{n}(\rho_{k-1})} 
    \end{bmatrix}\]
Where 
\[A = \begin{bmatrix}
    \frac{k-1}{k}\widehat{m_1}(\rho_0) & -\frac{1}{k}\widehat{m_1}(\rho_1) & ... & -\frac{1}{k}\widehat{m_1}(\rho_{k-1}) \\
    -\frac{1}{k}\widehat{m_1}(\rho_0) & \frac{k-1}{k}\widehat{m_1}(\rho_1) & .... & -\frac{1}{k}\widehat{m_1}(\rho_{k-1}) \\
    ... & ... & ... & ... \\
    -\frac{1}{k}\widehat{m_1}(\rho_0) & -\frac{1}{k}\widehat{m_1}(\rho_1) & .... & \frac{k-1}{k}\widehat{m_1}(\rho_{k-1}) \\
\end{bmatrix}\]
\[A = \begin{bmatrix}
    \frac{k-1}{k} & -\frac{1}{k} & ... & -\frac{1}{k} \\
    -\frac{1}{k} & \frac{k-1}{k} & .... & -\frac{1}{k} \\
    ... & ... & ... & ... \\
    -\frac{1}{k} & -\frac{1}{k} & .... & \frac{k-1}{k}\\
\end{bmatrix}\begin{bmatrix}
    \widehat{p}(\rho_0) & 0 & ... & 0 \\
    0 & \widehat{p}(\rho_1) & .... & 0 \\
    ... & ... & ... & ... \\
    0 & 0 & .... & \widehat{p}(\rho_{k-1})\\
\end{bmatrix}\]
Which means that we have that 
\[A^{n-1}\widehat{m_{1}} = \widehat{m_{n}}\]
Taking the powers of the above matrix is faster and it is generally easier to compute the eigenvalues of such a matrix. For large matrices, the above's eigenvalues and eigenvectors will faster to compute due to the fact that it is a diagonal matrix multiplied by a nilpotent matrix.
\section*{Spectral Methods}
Before we introduce the last part of this paper, we must discuss some basic lemmas. The first is that $G$ is $d$-regular, then $\lambda_1$ or the largest eigenvalue is equal to $d$. It's also true that $|\lambda_n| \leq d$ as well. 

Let $A$ be the adjaceny matrix of a graoh. Another lemma is that $\text{Trace}(A^k)$ counts the number of closed $k$ length closed walks in the graph. 
\[\sum_{i = 1}^V \lambda_i^k = \text{Trace}(A^k)\]

The next part of the paper will also only work on \textbf{Vertex-Transitive Graphs}. Vertex Transitive Graphs are graphs where there exists a graph automorphism $\phi$ such that $\phi(i) = j$ for any vertices i and j. In simple terms, it means the graph looks the same from every vertex. It also means that each node has the same amount of closed walks length k that start and end at that node. This quantity is $\frac{\text{Trace}(A^k)}{V}$ where $V$ is the number of vertices in the graph. \cite{Graph:2001} We aim to prove the following
\begin{theorem}
For simple random walks in vertex transitive graphs.
  \[\sum_{n = 1}^\infty P(\tau_{i,j} = n)t^n =  \frac{V \det_{j,i}(I - \frac{t}{d}A)}{\sum_{j = 1}^V \prod_{i \ne j} (1 - \frac{t}{d}\lambda_i)}\]
Where $A$ is the adjacency matrix of a such a graph $G$ and $d$ is the common degree of nodes in $G$. $\det_{ij}$ is the determinant of a matrix with the ith row and jth column removed.
\end{theorem}
 Let $M_{n_{i,j}} = P(\tau_{i,j} = n)$. 
It is clear that any vertex transitive graph that we have.
\[M_0 = I\]
and
\[M_1 = \frac{1}{d}A^1 - \frac{1}{d}\frac{Trace(A)}{V}I\]
This makes sure that a vertex can't visit itself in one move.
We also have.
\[M_2 = \frac{1}{d^2}A^2 - \frac{1}{d}\frac{Trace(A)}{V}M_1-\frac{1}{d^2}\frac{Trace(A^2)}{V}M_0\]
This has it so we can't visit our selves in 2 moves or 1 move.
Generalizing this
\[M_n = \frac{1}{d^n}A^n - \sum_{k = 1}^{n} \frac{1}{d^{k}} \frac{Trace(A^{k})}{V} M_{n - k} \]
Basically this says, we are not allowed to visit the node earlier than n moves and then make a loop to that same node to achieve a walk of technically \textit{technically} n moves from i to j. We get this sum.
\[\sum_{k = 0}^{n} \frac{1}{d^{k}} \frac{Trace(A^{k})}{V} M_{n -k} = \frac{1}{d^n}A^n  \]
As this is a Cauchy product, we get the following
\[\sum_{n = 0}^{\infty} \frac{t^n}{d^{n}} \frac{Trace(A^{n})}{V}\sum_{n=0}^\infty M_{n}t^n = \sum_{n = 0}^\infty \frac{t^n}{d^n}A^n  \]
By definition of trace we have.
\[\frac{1}{V}\left(\sum_{n = 0}^{\infty} \frac{t^k}{d^{k}}\sum_{i = 0}^V\lambda_i^k\right)\sum_{n=0}^\infty M_{n}t^n = \sum_{n = 0}^\infty \frac{t^n}{d^n}A^n  \]
\[\frac{1}{V}\left(\sum_{i = 0}^V\sum_{n = 0}^{\infty} \lambda_i^k\frac{t^k}{d^{k}}\right)\sum_{n=0}^\infty M_{n}t^n = \sum_{n = 0}^\infty \frac{t^n}{d^n}A^n  \]
\[\frac{1}{V}\left(\sum_{i = 0}^V \frac{1}{1 - \frac{\lambda_i t}{d}}\right)\sum_{n=0}^\infty M_{n}t^n = (I - \frac{t}{d}A)^{-1}  \]
so for $|t| < 1$
\[\sum_{n=0}^\infty M_{n}t^ n = \frac{V(I - \frac{t}{d}A)^{-1}}{\sum_{i = 0}^V \frac{1}{1 - \frac{\lambda_i t}{d}}}  \]

\[\sum_{n = 1}^\infty e_i^TM_nt^ne_j = \frac{V \frac{\det_{j,i}(I - \frac{t}{d}A)}{\prod_{i = 1}^V (1 - \frac{t}{d}\lambda_i)}}{\frac{\sum_{j = 1}^V \prod_{i \ne j} (1 - \frac{t}{d}\lambda_i)}{\prod_{i = 1}^V (1 - \frac{t}{d}\lambda_i)}}\]
\[\sum_{n = 1}^\infty P(\tau_{i,j} = n)t^n =  \frac{V \det_{j,i}(I - \frac{t}{d}A)}{\sum_{j = 1}^V \prod_{i \ne j} (1 - \frac{t}{d}\lambda_i)}\]
This is a very general result for all vertex transitive graphs which apply to all Cayley graphs as well.
\subsection*{Sample Distributions}
Using Mathematica, we can get series that describe our distributons using the general form above.
\begin{itemize}
    \item 3d Hypercube Walk from $i = (0,0,0) \to j = (1,1,1)$
    \[\sum_{n= 0}^\infty P(\tau_{i,j} = n)t^n =\frac{2 t^3}{9}+\frac{14 t^5}{81}+\frac{98 t^7}{729}+\frac{686 t^9}{6561}+\frac{4802 t^{11}}{59049}+O\left(t^{13}\right) \]
    \item Walk in $S_3$ from $i = e$ to $j = (1,3)$ 
    \[\sum_{n= 0}^\infty P(\tau_{i,j} = n)t^n = \frac{t}{3} + \frac{4 t^3}{27} + \frac{2 t^4}{27} + \frac{20 t^5}{243} + \frac{44 t^6}{729} + \frac{116 t^7}{2187} + \frac{280 t^8}{6561} + O(t^{9})\]
    \item Walk in $D_8$ from $i = e$ to $j = (14)(23)$
    \[\sum_{n= 0}^\infty P(\tau_{i,j} = n)t^n = \frac{t}{3}+\frac{4 t^3}{27}+\frac{28 t^5}{243}+\frac{196 t^7}{2187}+\frac{1372 t^9}{19683}+\frac{9604 t^{11}}{177147}+O\left(t^{13}\right)\]
    It is important to note that the generators of $S_3$ in this example are the transpositions and that generators of $D_8$ are $(14)(23)$ and $(1234)$ and $(1 2)(3 4)$. As above, we walk along our generators with equal probability. 
\end{itemize}
\section*{Continuous-Time Random Walks}

Lastly, we consider a continuous-time random walk on a finite state space, where transitions between states occur according to a Poisson process with rate \( \lambda = 1 \). Let $M$ be the markov matrix defining this walk. Let $j$ be the defined as the end point. Let $Q$ be the substochastic matrix of $M$ with the $j$th row and column removed. Let $P_1$ again be the vector with the 1 step probabilities.

Define \( \tau^c_{*,j} \) as the vector random variable representing the first hitting times to node \( j \) from every other node. We are interested in computing the cumulative distribution function (CDF) of \( \tau^c_{*,j} \), that is,
\[
P(\tau^c_{*,j} \leq t).
\]

In continuous time, the probability that the walk makes exactly \( n \) transitions by time \( t \) is given by the Poisson distribution:
\[
P(\text{makes } n \text{ steps}) = \frac{t^n e^{-t}}{n!}.
\]

Therefore, the probability that the walk hits node \( j \) within time \( t \) can be written as
\[
P(\tau^c_{*,j} \leq t) = \sum_{n=0}^\infty \frac{t^n e^{-t}}{n!} \cdot P(\text{hits } j \text{ in at most } n \text{ steps}).
\]

The second term in the sum can be written as
\[
\sum_{k = 1}^{n} Q^{k-1} P_1,
\]
Hence,
\[
P(\tau^c_{*,j} \leq t) = \sum_{n=0}^\infty \frac{t^n e^{-t}}{n!} \sum_{k=1}^n Q^{k-1} P_1.
\]

We can interchange the order of summation:
\[
P(\tau^c_{*,j} \leq t) = \sum_{k=1}^\infty Q^{k-1} P_1 \sum_{n = k}^\infty \frac{t^n e^{-t}}{n!}.
\]

Recognizing that this is the tail of the Poisson distribution, we move to a matrix formulation. Using the identity
\[
\sum_{n=0}^\infty \frac{t^n e^{-t}}{n!} Q^{n} = e^{-t(I - Q)},
\]
we obtain a closed-form expression:
\[
P(\tau^c_{*,j} \leq t) = (I - Q)^{-1} \left( I - e^{-t(I - Q)} \right) P_1.
\text{ for } t > 0\]

Differentiating this expression with respect to \( t \), we obtain the vector probability density function (PDF):
\[
f_{\tau^c_{*,j}}(t) = \frac{d}{dt} P(\tau^c_{*,j} \leq t) = e^{-t(I - Q)} P_1.
\]
We can interpet this as smoothed out version of our orginial walk. This can make the moments easier to obtain as the above distribution has the same moments as our discrete time random walk.
\[\int_{0}^\infty te^{-t(I - Q)} P_1 = (I-Q)^{-2}P_1 = \sum_{n = 0}^\infty nQ^{n-1}P_1\]
\[\int_{0}^\infty t^2e^{-t(I - Q)} P_1 = 2(I-Q)^{-3}P_1 = \sum_{n = 0}^\infty n^2Q^{n-1}P_1\]
And so on. 
\section*{Future Work}
The main future work is obtaining distributions from $M_n$ by extracting coefficients through analysis methods.

Another direction one could take is by calculating how much information one would need to have about the $M_n$ matrices to get the adjacency matrix of the graph. It is clear that recovering the matrix from the distribution of one pair of vertices is not sufficient (consider cycle and path graphs). This begs the question of how many distributions are needed. Since we have that $A^k$ are in the Bose-Mesner Algebra by the following recursion explored earlier, 
\[\sum_{k = 0}^{n} \frac{1}{d^{k}} \frac{Trace(A^{k})}{V} M_{n -k} = \frac{1}{d^n}A^n  \]
we can retrieve $A$ from $M_n$. Then, therefore there exists some proportion of $M_n$ that is needed before extracting $A$.
\section*{Acknowledgments}
Thank you to Jonathan Pakianthan, Alex Iosevich and Matthew Dannenburg for their continued support.  Thank you to my team members at STEM FOR ALL 2024: Yujia Hu, Yiling Zou, Charlie Li. Thank you to my advisor Sushil Varma.
\bibliographystyle{plain} 
\bibliography{references} 
\pagebreak
\section*{Appendix}
\subsection*{Proof of Theorem 2}
\[\phi_{\tau_{*,j}}(t) = e^{it}P_1\]
\[(I - e^{it}Q)\phi_{\tau_{*,j}}(t)= e^{it}P_1\]
We take the derivative of both sides
\[(I - e^{it}Q)\phi'_{\tau_{*,j}}(t) - ie^{it}Q\phi_{\tau_{*,j}}(t) = ie^{it}P_1\]
\[(I - e^{it}Q) \phi'_{\tau_{*,j}}(t)= ie^{it}P_1 + ie^{it}Q\phi_{\tau_{*,j}}(t)\]
We now plug in $t =0$. 
\[(I - Q) \phi'_{\tau_{*,j}}(0)= iP_1 + iQ\phi_{\tau_{*,j}}(0)\]
Since $\phi_{\tau_{i,j}}(0) =  \sum_{n = 1}^\infty e^{i0t}P(\tau_{i,j} = n) = 1$. If $\mathbf{1}$ is the vector of all 1's, then we have that
\[(I - Q)\phi'_{\tau_{*,j}}(0) = iP_1 + Q\mathbf{1}\]
\[\phi'_{\tau_{*,j}}(0) = i(I - Q)^{-1}(P_1 + Q\mathbf{1})\]
It follows that $Q\mathbf{1}_i = \sum_{k \ne j}P(\tau_{i,k} = 1)$. Therefore, $(P_1 + Q\mathbf{1}) = \mathbf{1}$
\[\phi'_{\tau_{*,j}}(0) = i(I - Q)^{-1}\mathbf{1}\]
So then
\[E[\tau_{*,j}] = (I - Q)^{-1}\mathbf{1} = \sum_{n = 0}^{\infty} Q^n \mathbf{1}\]
We then have that second derivative of this is
\[\phi''_{\tau_{*,j}}(t)= -2Q(I-Q)^{-2}\mathbf{1}\]
Which means that 
\[E[\tau_{*,j}^2] = 2Q(I-Q)^{-2}\mathbf{1} \]
\[E[\tau_{*,j}^2] = 2\sum_{n=0}^\infty nQ^n\mathbf{1} \]
\subsection*{The Distribution of the Cycle}
\subsection*{Cycle Case} Let's perform these computations for a cycle. It follows for a k-cycle (denoted as $C_k$) \[Q = \begin{bmatrix} 0 & \frac{1}{2} & 0 & 0 & ... & 0 & 0 \\ \frac{1}{2} & 0 & \frac{1}{2} & 0 &... &0 &0 \\ 0 & \frac{1}{2} & 0 & \frac{1}{2} &... &0 &0 \\ ... & ... & ... & ... & ... &... &... \\ 0 & 0 & 0 & 0 & ... & \frac{1}{2} & 0 \\ \end{bmatrix},P_1 = \begin{bmatrix} \frac{1}{2}\\0 \\\vdots \\ 0 \\ \frac{1}{2} \end{bmatrix}\] So then we have \[\begin{bmatrix} 0 & \frac{1}{2} & 0 & 0 & ... & 0 & 0 \\ \frac{1}{2} & 0 & \frac{1}{2} & 0 &... &0 &0 \\ 0 & \frac{1}{2} & 0 & \frac{1}{2} &... &0 &0 \\ ... & ... & ... & ... & ... &... &... \\ 0 & 0 & 0 & 0 & ... & \frac{1}{2} & 0 \\ \end{bmatrix}\begin{bmatrix} P(\tau_{1,0} = n-1)\\ P(\tau_{2,0} = n-1) \\ ... \\ P(\tau_{k-1,0} = n-1) \end{bmatrix} = \begin{bmatrix} P(\tau_{1,0} = n)\\ P(\tau_{2,0} = n) \\ ... \\ P(\tau_{k-1,0} = n) \end{bmatrix}\] Or in other words \[\begin{bmatrix} 0 & \frac{1}{2} & 0 & 0 & ... & 0 & 0 \\ \frac{1}{2} & 0 & \frac{1}{2} & 0 &... &0 &0 \\ 0 & \frac{1}{2} & 0 & \frac{1}{2} &... &0 &0 \\ ... & ... & ... & ... & ... &... &... \\ 0 & 0 & 0 & 0 & ... & \frac{1}{2} & 0 \\ \end{bmatrix}^{n-1}\begin{bmatrix} P(\tau_{1,0} = 1)\\ P(\tau_{2,0} = 1) \\ ... \\ P(\tau_{k-1,0} = 1) \end{bmatrix} = \begin{bmatrix} m_{n}(1) \\ m_{n}(2) \\ ... \\ m_{n}(k-1) \end{bmatrix}\] So it follows we want to take the diagonalization of the Toeplitz Matrix above. Let us call that matrix $H$ and its diagonlization $LDL^{-1}$. Thankfully, the eigenvectors and eigenvalues for tridiagonal toeplitz matrices are well known and with that we have the following as the diagonalization. These eigenvalues were found in \cite{duru2018toeplitz} \[\begin{bmatrix} 0 & \frac{1}{2} & 0 & 0 & ... & 0 & 0 \\ \frac{1}{2} & 0 & \frac{1}{2} & 0 &... &0 &0 \\ 0 & \frac{1}{2} & 0 & \frac{1}{2} &... &0 &0 \\ ... & ... & ... & ... & ... &... &... \\ 0 & 0 & 0 & 0 & ... & \frac{1}{2} & 0 \\ \end{bmatrix}\] \[= \frac{2}{k}\begin{bmatrix} \sin(\frac{\pi}{k})& \sin(\frac{2\pi}{k}) & ... & \sin(\frac{(k-1)\pi}{k}) \\ \sin(\frac{2\pi}{k}) & \sin(\frac{4\pi}{k}) & ... & \sin(\frac{2(k-1)\pi}{k}) \\ ... & ... & ... & ... \\ \sin(\frac{(k-1)\pi}{k})& \sin(\frac{2(k-1)\pi}{k})& ... & \sin(\frac{(k-1)^2\pi}{k}) \\ \end{bmatrix}\begin{bmatrix} \cos(\frac{\pi}{k})& & 0 & ... & 0 \\ 0 & \cos(\frac{2\pi}{k}) & 0 & ... &0 \\ ... & ... & ... & ... & ... \\ 0 & 0 & 0 & ... & \cos(\frac{(k-1)\pi}{k}) \\ \end{bmatrix}\]\[\begin{bmatrix} \sin(\frac{\pi}{k})& \sin(\frac{2\pi}{k}) & ... & \sin(\frac{(k-1)\pi}{k}) \\ \sin(\frac{2\pi}{k}) & \sin(\frac{4\pi}{k}) & ... & \sin(\frac{2(k-1)\pi}{k}) \\ ... & ... & ... & ... \\ \sin(\frac{(k-1)\pi}{k})& \sin(\frac{2(k-1)\pi}{k})& ... & \sin(\frac{(k-1)^2\pi}{k}) \\ \end{bmatrix}\] As we know from earlier, \[\begin{bmatrix} P(\tau_{1,0} = 1)\\ P(\tau_{2,0} = 1) \\ ... \\ P(\tau_{k-1,0} = 1) \end{bmatrix}= \begin{bmatrix} \frac{1}{2} \\ 0 \\ ... \\ \frac{1}{2} \end{bmatrix}\] So, then we first multiply by $L^{-1}$ or the last matrix in the diagonalization. \[\begin{bmatrix} \sin(\frac{\pi}{k})& \sin(\frac{2\pi}{k}) & ... & \sin(\frac{(k-1)\pi}{k}) \\ \sin(\frac{2\pi}{k}) & \sin(\frac{4\pi}{k}) & ... & \sin(\frac{2(k-1)\pi}{k}) \\ ... & ... & ... & ... \\ \sin(\frac{(k-1)\pi}{k})& \sin(\frac{2(k-1)\pi}{k})& ... & \sin(\frac{(k-1)^2\pi}{k}) \\ \end{bmatrix}\begin{bmatrix}\frac{1}{2} \\ 0 \\ ... \\ \frac{1}{2} \end{bmatrix} = \frac{1}{2} \begin{bmatrix} \sin(\frac{\pi}{k}) + \sin(\frac{(k-1)\pi}{k})\\ \sin(\frac{2\pi}{k}) + \sin(\frac{2(k-1)\pi}{k}) \\ ... \\ \sin(\frac{(k-1)\pi}{k}) + \sin(\frac{(k-1)^2\pi}{k}) \end{bmatrix}\] We next multiply by the diagonal matrix. \[\frac{1}{2}\begin{bmatrix} \cos(\frac{\pi}{k})& & 0 & ... & 0 \\ 0 & \cos(\frac{2\pi}{k}) & 0 & ... &0 \\ ... & ... & ... & ... & ... \\ 0 & 0 & 0 & ... & \cos(\frac{(k-1)\pi}{k}) \\ \end{bmatrix}^{n-1}\begin{bmatrix} 0\\ \sin(\frac{2\pi}{k}) \\ 0 \\ \sin(\frac{4\pi}{k}) \\ ... \\ \end{bmatrix} =\frac{1}{2} \begin{bmatrix} cos(\frac{\pi}{k})\sin(\frac{\pi}{k})^{n-1} + \sin(\frac{(k-1)\pi}{k}))\\ cos(\frac{2\pi}{k})^{n-1}\sin(\frac{2\pi}{k}) + \sin(\frac{2(k-1)\pi}{k}) \\ ... \\ cos(\frac{(k-1)\pi}{k})^{n-1}(\sin(\frac{(k-1)\pi}{k}) + \sin(\frac{(k-1)^2\pi}{k}) ) \end{bmatrix}\] And lastly multiplying by $L$. \[\frac{2}{k}\begin{bmatrix} \sin(\frac{\pi}{k})& \sin(\frac{2\pi}{k}) & ... & \sin(\frac{(k-1)\pi}{k}) \\ \sin(\frac{2\pi}{k}) & \sin(\frac{4\pi}{k}) & ... & \sin(\frac{2(k-1)\pi}{k}) \\ ... & ... & ... & ... \\ \sin(\frac{(k-1)\pi}{k})& \sin(\frac{2(k-1)\pi}{k})& ... & \sin(\frac{(k-1)^2\pi}{k}) \\ \end{bmatrix}\frac{1}{2} \begin{bmatrix} cos(\frac{\pi}{k})\sin(\frac{\pi}{k})^{n-1} + \sin(\frac{(k-1)\pi}{k}))\\ cos(\frac{2\pi}{k})^{n-1}\sin(\frac{2\pi}{k}) + \sin(\frac{2(k-1)\pi}{k}) \\ ... \\ cos(\frac{(k-1)\pi}{k})^{n-1}(\sin(\frac{(k-1)\pi}{k}) + \sin(\frac{(k-1)^2\pi}{k}) ) \end{bmatrix}\] \[= \begin{bmatrix} \frac{1}{k}\sum_{j = 0}^{k-1} \cos(\frac{j\pi}{k})^{n-1}(\sin(\frac{j\pi}{k}) + \sin(\frac{j(k-1)\pi}{k}))\sin(\frac{j\pi}{k}) \\ \frac{1}{k}\sum_{j = 0}^{k-1} \cos(\frac{j\pi}{k})^{n-1}(\sin(\frac{j\pi}{k}) + \sin(\frac{j(k-1)\pi}{k}))\sin(\frac{2j\pi}{k}) \\ \frac{1}{k}\sum_{j = 0}^{k-1} \cos(\frac{j\pi}{k})^{n-1}(\sin(\frac{j\pi}{k}) + \sin(\frac{j(k-1)\pi}{k}))\sin(\frac{3j\pi}{k}) \\ ... \\ \frac{1}{k}\sum_{j = 0}^{k-1} \cos(\frac{j\pi}{k})^{n-1}(\sin(\frac{j\pi}{k}) + \sin(\frac{j(k-1)\pi}{k}))\sin(\frac{(k-1)j\pi}{k}) \end{bmatrix}\] Therefore \[P(\tau_{i,j} = n) = \frac{1}{k}\sum_{j = 0}^{k-1} \cos(\frac{j\pi}{k})^{n-1}(\sin(\frac{j\pi}{k}) + \sin(\frac{j(k-1)\pi}{k}))\sin(\frac{ij\pi}{k}) \]
\end{document}